\documentclass[reprint,aps,prx,floatfix,superscriptaddress]{revtex4-2-shared-email}
\usepackage{amsmath}
\usepackage{amsfonts}
\usepackage{physics}
\usepackage{graphicx}
\usepackage[breaklinks = true,
            colorlinks = true,
            linkcolor  = blue,
            citecolor  = blue,
            urlcolor   = blue]
            {hyperref}

\usepackage[capitalise,nameinlink]{cleveref}
\crefname{figure}{Fig.}{Figs.}
\Crefname{figure}{Figure}{Figures}
\crefname{table}{Table}{Tables}
\Crefname{table}{Table}{Tables}
\crefname{section}{Sec.}{Secs.}
\Crefname{section}{Section}{Sections}
\crefname{subsection}{Sec.}{Secs.}
\Crefname{subsection}{Section}{Sections}
\crefname{equation}{Eq.}{Eqs.}
\Crefname{equation}{Equation}{Equations}
\crefname{appendix}{Appendix}{Appendices}
\Crefname{appendix}{Appendix}{Appendices}

\usepackage{etoolbox}
\apptocmd{\sloppy}{\hbadness 10000\relax}{}{}

\begin{document}
\title{Tensor-Network Finite Elements for Analytic Operator Equations} 

\author{Abhijatmedhi Chotrattanapituk}
\thanks{These authors contributed equally}
\email{c\_earth@mit.edu}
\affiliation{Quantum Measurement Group, MIT, Cambridge, Massachusetts 02139, USA}
\affiliation{Department of Electrical Engineering and Computer Science, MIT, Cambridge, Massachusetts 02139, USA}

\author{Michael J. Landry}
\thanks{These authors contributed equally}
\affiliation{Quantum Measurement Group, MIT, Cambridge, Massachusetts 02139, USA}
\affiliation{Department of Nuclear Science and Engineering, MIT, Cambridge, Massachusetts 02139, USA}
\affiliation{Department of Physics, MIT, Cambridge, MA 02139, USA}

\author{Chu-Liang Fu}
\affiliation{Quantum Measurement Group, MIT, Cambridge, Massachusetts 02139, USA}
\affiliation{Department of Nuclear Science and Engineering, MIT, Cambridge, Massachusetts 02139, USA}

\author{Mingda Li}
\email{mingda@mit.edu}
\affiliation{Quantum Measurement Group, MIT, Cambridge, Massachusetts 02139, USA}
\affiliation{Department of Nuclear Science and Engineering, MIT, Cambridge, Massachusetts 02139, USA}

\begin{abstract}
Operator equations (OEs) underpin quantitative modeling across science and engineering. Finite-element (FE) methods discretize continuous OEs into finite-dimensional algebraic systems, whereas tensor networks (TNs) provide flexible variational representations of correlated discrete systems. Here, we develop a framework that connects FE with TN for analytic OEs. The power of this method comes from its ability to convert highly non-linear partial differential equations into linear matrix equations. In particular, we show that FE discretization induces a hierarchy of multilinear interaction tensors, through which differential, integral, nonlinear, memory, and delay equations can be expressed within a common algebraic structure. The resulting systems are reformulated as weighted-residual optimization problems over TN degrees of freedom. Matrix-product-state calculations for one-dimensional linear and nonlinear diffusion reproduce conventional solutions with controlled error while preserving continuity and Neumann boundary conditions. The framework provides a common variational language for analytic OEs and establishes a direct connection between FE numerical formalism and TN variational algorithms, offering a general foundation for TN-based and quantum-inspired approaches to solving OEs.
\end{abstract}

\maketitle

\section{Introduction}\label{sec:intro}
Operator equations (OEs) provide a general mathematical framework for quantitative modeling across science and engineering, which includes partial differential equations (PDEs), integro-differential equations (IDEs), and other nonlocal evolution quations~\cite{krasnoselskii1972approximate, atkinson2009theoretical, brezis2011functional}. Such equations arise naturally in electrodynamics, fluid dynamics, transport, quantum mechanics, statistical physics, and multi-scale materials modeling. However, only limited classes of OEs admit analytical solutions; consequently, numerical approximation methods have become indispensable for practical applications. Established discretization frameworks, including finite-difference (FD)~\cite{leveque2007finite, strikwerda2004finite, thomas1995numerical}, finite-element (FE)~\cite{ciarlet2002finite, brenner2008mathematical, hughes2000finite, zienkiewicz2013finite}, and spectral-element (SE) methods, can achieve high numerical accuracy, but their computational cost grows substantially as the discretization is refined and the number of degrees of freedom increases, particularly for large-scale, high-dimensional, and highly correlated systems. Among these approaches, the finite-element method (FEM) has become one of the most widely used discretization frameworks for continuous OEs~\cite{ciarlet2002finite, brenner2008mathematical}. In conventional FEM, continuous OEs are reduced to finite-dimensional algebraic systems by expanding the solution in localized basis functions and projecting the governing equation onto the associated approximation space through a weak formulation.

Tensor networks (TNs) provide efficient variational representations of high-dimensional tensors with nontrivial correlation structure through interconnected low-rank decompositions~\cite{orus_2014}. Originally developed in quantum many-body physics, TN models such as matrix product states (MPS)~\cite{Fannes_1992,Klumper_1991,Klumper_1993}, projected entangled pair states (PEPS)~\cite{verstraete_2004_A}, and multiscale entanglement renormalization ansatz (MERA)~\cite{Vidal_2007_B}, together with associated optimization techniques including the density matrix renormalization group (DMRG)~\cite{white_1992,white_1993}, have become powerful tools for representing and optimizing correlated high-dimensional systems. A central feature of TN methods is their ability to capture dominant correlation structures efficiently, often avoiding the full exponential complexity of explicit tensor representations for structured low-correlation systems. Owing to these capabilities, TN methods have found broad applications across quantum many-body physics~\cite{Jiang_2008,Picot_2016,Corboz_2016,Dubail_2015,Aguado_2008,Buerschaper_2009,Gu_2009}, quantum information~\cite{Markov_2008,Ferris_2014,Huggins_2019}, machine learning~\cite{Stoudenmire_2016,Han_2018_B,Levine_2019,Liu_2019,Gallego_2022}, quantum chemistry~\cite{Chan_2008,Szalay_2015,Krumnow_2016}, disordered systems~\cite{Paredes_2005,Chandran_2015,Pollmann_2016}, and quantum field theory~\cite{Verstraete_2010,Haegeman_2013_A,Swingle_2012}. Recent works have also explored TN approaches for PDE solving~\cite{Bachmayr_2016,Gourianov_2022,Ye_2022}, where TN representations are primarily employed to compress discretized solution spaces and reduce computational cost~\cite{Lubasch_2018,Garca_Ripoll_2021,Cao_2013,Costa_2019,Ameri_2023,Liu_2021,Kyriienko_2021}. 

Building on these developments, we develop a different approach to solving OEs with TNs. Rather than treating TNs as compressed surrogates of numerical solutions, we directly formulate discretized OEs in terms of TN degrees of freedom. Specifically, we generalize the conventional FEM formulation into a tensor-product representation. Through this construction, the conventional FE coefficient tensor is replaced by a structured coefficient TN. The discretized OE is then reformulated as an optimization problem over TN parameters. In this sense, TNs serve not merely as compression tools, but as variational representations of correlated FE coefficient states. Our formalism offers three key advantages: it provides a unified treatment of different classes of OEs, incorporates inter-element correlations directly into the variational space, and enables established TN algorithms to operate directly on discretized equations. More broadly, it establishes a direct connection between FEM numerical methods and TN variational algorithms, providing a general foundation for TN-based and quantum-inspired approaches to solving OEs.

The remainder of this paper is organized as follows. In \cref{sec:note}, we introduce the mathematical notations used throughout this work. In \cref{sec:fe}, we formulate analytic OEs and derive their FE discretization. In \cref{sec:tn}, we develop the TN formulation from FE discretization of analytic OEs. In \cref{sec:reduc}, we show how the general formalism specializes to some important sub-classes of OEs. In \cref{sec:demo}, we present representative demonstrations of the method. Finally, in \cref{sec:conc}, we discuss implications, limitations, and possible extensions of the proposed approach. Additional derivations and numerical results are provided in the Appendices.
\section{Notations}\label{sec:note}
Non-scalar tensors are denoted by bold Roman characters, while their components (including scalar tensors) are written in italic font with Greek tensor indices. For example, for a rank-two contravariant tensor $\mathbf a$, its $(\mu,\nu)$ component is denoted by $a^{\mu\nu}$. Non-tensor indices are denoted by upright Roman characters. Einstein summation convention is also assumed over repeated tensor indices.

We use bracket notation to denote sequences generated by all values of the free indices appearing inside the brackets. For any tensor component, bracket notation simply recover the tensor, e.g, $\left[a^{\mu\nu}\right]=\mathbf{a}$. For tensors with many indices, bracket notation may also denote index strings compactly as $a^{\alpha_1\alpha_2\cdots} = a^{[\alpha_\mathrm{i}]}$ where the sequence notation does not apply to $\alpha$ by context. For situations in which we wish to be explicit, we will use notation of the form
\begin{equation*}
    \partial_{\mathbf{b}}\left(a^{[\alpha_\mathrm{i}]}\right)
    = \left[\frac{\partial}{\partial b^{[\beta_\mathrm{j}]}} a^{[\alpha_\mathrm{i}]}\right]_{[\beta_\mathrm{j}]}
\end{equation*}
where $\mathbf{a} =\left[a^{[\alpha_\mathrm{i}]}\right]$, $\mathbf{b}=\left[b^{[\beta_\mathrm{j}]}\right]$.

We also reserve ``$\otimes$'' for tensor product (TP), ``$\odot$'' for tensor contraction, ``$\oplus$'' for direct sum.
\section{Finite-Element Discretization of Operator Equations}\label{sec:fe}
FEM reduces continuous OEs to finite-dimensional algebraic systems by approximating the solution within a finite basis and projecting the governing equation onto the approximation space through a weak formulation, which induces a hierarchy of multi-linear interaction tensors. These tensors encode the interaction geometry of the original OE, including locality, memory effects, stochastic couplings, and higher-order correlations. It therefore provides the natural starting point for our construction, converting the continuous OE into finite-dimensional residual tensors over the FE coefficients that can subsequently be represented and optimized within a TN variational space.

\subsection{Weak Formulation}
Let $\Omega$ denote the domain of the independent variable $\mathbf{q}$, $\mathcal{H}(\Omega)$ denote the corresponding admissible solution space, and $\mathcal{R}(\Omega)$ denote the residual space, whose elements quantify the extent to which a trial solution fails to satisfy the governing equation. We write a general OE over $\Omega$ as 
\begin{equation}
    \mathbf{L}[\mathbf{u}] = \mathbf{0}
\end{equation}
where $\mathbf{L}:\mathcal{H}(\Omega)\rightarrow\mathcal{R}(\Omega)$, and $\mathbf{u}\in\mathcal{H}(\Omega)$. Then, consider a finite collection of scalar basis functions $b_\beta:\Omega\rightarrow \mathbb{F}$ where $\mathbb{F}$ indicates the scalar field. In FE representation, the solution can be approximated as
\begin{equation}\label{eq:sol_approx}
    \mathbf{u} \approx \mathbf{c}^\beta b_\beta
\end{equation}
where $\left[\mathbf{c}^\beta\right]$ is the coefficient tensor. The Galerkin weak formulation would requires
\begin{equation}
    \int_\Omega b^\ast_\beta\mathbf{L}[\mathbf{u}]\, \mathrm{d}\mathbf{q}
    = \mathbf{0}
\end{equation}
for every test function $b_\beta$ with decorator $\ast$ indicates conjugation. This reduces the continuous equation to a finite system of algebraic equations for the coefficient tensor.

\subsection{Analytic Operator Equations}
\begin{figure}[!htb]
    \centering
    \includegraphics[width=0.9\columnwidth]{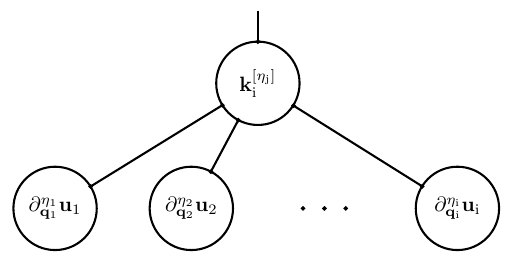}
    \caption{Contraction of $\mathrm{i}$-linear interaction kernel $\mathbf{k}_{\mathrm{i}}^{\left[\eta_\mathrm{j}\right]}$ with an $\mathrm{i}$-point correlation in tensor network representation. The solid objects represent tensors, hanging lines represent free indices, lines connecting two tensors represents index contractions. The rank of the kernel is chosen such that it can contract each of the derivative term, $\partial_{\mathbf{q}_\mathrm{j}}^{\eta_\mathrm{j}}\mathbf{u}_\mathrm{j}$, completely with the left over rank (top branch) the same as $\mathbf{L}[\mathbf{u}]$.}
    \label{fig:interact_kernel}
\end{figure}

We consider a general class of analytic OEs that admit a convergent Fréchet–Taylor expansion~\cite{boyd1984analytical,boyd1985fading,schetzen1980volterra,mujica1986complex,dineen1999complex} of the form
\begin{equation}\label{eq:aoe}
    \mathbf{L}[\mathbf{u}]
    = \sum_{\mathrm{i}=0}^\infty\dfrac{1}{\mathrm{i}!}\int_{\Omega^\mathrm{i}}\mathbf{k}_{\mathrm{i}}^{\left[\eta_\mathrm{j}\right]}\odot\bigotimes_{\mathrm{j}=1}^\mathrm{i}\partial_{\mathbf{q}_\mathrm{j}}^{\eta_\mathrm{j}}\mathbf{u}_\mathrm{j}\,\mathrm{d}\mathbf{q}_\mathrm{j}
\end{equation}
where $\eta_\mathrm{j}\in\mathbb{N}_0$, $\mathbf{u}_\mathrm{j} = \mathbf{u}(\mathbf{q}_\mathrm{j})$, and $\mathbf{k}_{\mathrm{i}}^{\left[\eta_\mathrm{j}\right]}$ is an $\mathrm{i}$-linear interaction kernel with proper contraction to each $\partial_{\mathbf{q}_\mathrm{j}}^{\eta_\mathrm{j}}\mathbf{u}_\mathrm{j}$ as shown in \cref{fig:interact_kernel}. 

If we further define
\begin{equation}\label{eq:con_oe}
    \mathbf{f}_{\mathrm{i}, \left[\beta_\mathrm{j}\right]}
    = \dfrac{1}{\mathrm{i}!}\int_{\Omega^\mathrm{i}}\mathbf{k}_{\mathrm{i}}^{\left[\eta_\mathrm{j}\right]}\odot\bigotimes_{\mathrm{j}=1}^\mathrm{i}\partial_{\mathbf{q}_\mathrm{j}}^{\eta_\mathrm{j}}b_{\beta_\mathrm{j}}\,\mathrm{d}\mathbf{q}_\mathrm{j}
\end{equation}
where $b_{\beta_\mathrm{j}}=b_{\beta_\mathrm{j}}(\mathbf{q}_\mathrm{j})$, the weak formulation can be written, for all $\beta$, as
\begin{equation}\label{eq:wf_aoe}
    \sum_{\mathrm{i}=0}^\infty \mathbf{g}_{\mathrm{i},\beta, \left[\beta_\mathrm{j}\right]}\odot\bigotimes_{\mathrm{j}=1}^\mathrm{i}\mathbf{c}^{\beta_\mathrm{j}}
    = \mathbf{0}
\end{equation}
where the discretized operator tensor is defined as
\begin{equation}
    \mathbf{g}_{\mathrm{i},\beta,\left[\beta_\mathrm{j}\right]} = \int_\Omega b^\ast_\beta\mathbf{f}_{\mathrm{i}, \left[\beta_\mathrm{j}\right]}\, \mathrm{d}\mathbf{q}\,.
\end{equation}

\subsection{Finite-Element}
Suppose that $\Omega$ is partitioned into $N$ finite elements $\Omega_\mathrm{n}$, where $\mathrm{n}\in\{1,\ldots,N\}$, each with local basis functions $b_{\mathrm{n},\beta_\mathrm{n}}$ which are non-zero only in $\Omega_\mathrm{n}$, and elements can have different numbers of basis functions. For cleanliness, we will slightly abuse our notation by acting as if all elements contain the same number of basis functions. With this assumption, we can substitute $\mathrm{n}$ with a tensor index $\nu$, and we can replace \cref{eq:sol_approx} with
\begin{equation}
    \mathbf{u} \approx \mathbf{c}^{\nu\beta} b_{\nu\beta}\,.
\end{equation}
We can cascade this to the case for analytic OEs, and get the replacement for \cref{eq:wf_aoe} as
\begin{equation}\label{eq:wf_fe_aoe}
    \sum_{\mathrm{i}} \mathbf{g}_{\mathrm{i},\nu\beta, \left[\nu_\mathrm{j}\beta_\mathrm{j}\right]}\odot\bigotimes_{\mathrm{j}=1}^\mathrm{i}\mathbf{c}^{\nu_\mathrm{j}\beta_\mathrm{j}}
    = \mathbf{0}
\end{equation}
where $\mathbf{g}_{\mathrm{i},\nu\beta, \left[\nu_\mathrm{j}\beta_\mathrm{j}\right]}$ can be obtain from simple substitution of basis function indexings.
\section{Tensor Network Formulation of Discretized Analytic Operator Equations}\label{sec:tn}

Our goal here is to convert a non-linear PDE into a linear matrix equation. To see what ingredients will be needed, begin by considering a much simpler problem, namely converting a linear PDE into a matrix equation. The solution is essentially trivial: any function can be represented as a linear combination of basis elements in a Hilbert space of square-integrable functions. (In practice, we will truncate the infinite set of basis elements and define the \emph{approximation space} to be the span of the retained basis elements.) In this basis, the differential operator of the PDE can be readily converted into a matrix. If the PDE is nonlinear, however, we must be more careful. There is no clean Hilbert space representation. Instead, we must upgrade the Hilbert space to Fock space. In this way each state-vector in Fock space contains multiple copies of the Hilbert space vector, meaning that operations like multiplication of the function by e.g. itself or its derivative can then be represented as linear operators acting on Fock space, subject to certain constraints. The result is a set of constrained linear equations, which are substantially easier to solve than the original non-linear PDE. 

The procedure for mapping non-linear PDEs to linear matrix equations on Fock space proceeds as follows. The FE discretization reduces an analytic OE to a finite-dimensional tensor system for the coefficient tensor. In conventional FEM, the solution space is represented as direct sum of all local basis functions, $b_{\mathrm{n},\beta_\mathrm{n}}$, which are used to construct the Hilbert space and, subsequently, the Fock space. As each of them is non-zero in only one element, the coefficients associated with different elements are independent variational degrees of freedom, and correlations between elements enter indirectly through the discretized operator tensors. In this section, we enlarge the coefficient representation so that correlations between finite elements can be encoded directly.

\subsection{Augmented Element Spaces}

The local approximation space on element $\mathrm{n}$ is the span of all local basis functions of that element. For normal FEM, the overall approximation Hilbert space $\Phi\subset\mathcal{H}(\Omega)$ is the direct sum of the local approximation spaces:
\begin{equation}
    \Phi = \bigoplus_{\mathrm{n}=1}^N\Phi_\mathrm{n}
\end{equation}
where 
\begin{equation}\label{eq:local_approx_space}
    \Phi_\mathrm{n} = \text{Span}\left(\left\{b_{\mathrm{n},\beta_\mathrm{n}}\mid \beta_\mathrm{n}\in\{1,\ldots,B_\mathrm{n}\}\right\}\right)
\end{equation}
and $B_\mathrm{n}$ is the number of basis functions in element $\mathrm{n}$. To build a TP representation, we use bra-ket notation from quantum mechanics, replacing each basis element by
\begin{equation}
b_{\mathrm{n},\beta_\mathrm{n}} \mapsto \ket{b_{\mathrm{n},\beta_\mathrm{n}}} \in \Phi\,.
\end{equation}
If we were only interested in linear PDEs, we could stop here. All such equations can (trivially) be projected down to matrix equations acting on elements of the approximation Hilbert space $\Phi$. We are, however, primarily interested in the case of nonlinear PDEs. As such we define the Fock space by 
\begin{equation}\label{eq:augmented_fe_space}
    \Psi
    = \bigoplus_{\mathrm{i}=0}^\infty \Phi^{\otimes\mathrm{i}}\,.
\end{equation}
This space has a choice of basis vectors of the form
\begin{equation}
    \ket{\left[\mathrm{n}_\mathrm{j}\beta_{\mathrm{n}_\mathrm{j}}\right]}_\mathrm{i}=\ket{\mathrm{n}_1\beta_{\mathrm{n}_1},\ldots,\mathrm{n}_\mathrm{i}\beta_{\mathrm{n}_\mathrm{i}}}\in\Psi\,.
\end{equation}

By itself, the Fock space is not identified directly with approximation space over the physical domain. It serves as a space for variational coefficients. To obtain physical solution, one must define a projection. For the purpose of solving OEs, we will choose the projection from our physical interpretation of states in $\Psi$ where each basis vector represents a sequence of basis functions involve in a correlation term, e.g., if $\mathbf{u} \approx \widetilde{\mathbf{c}}^{\nu\beta} b_{\nu\beta}$, we will replace the coefficient for the standard 2-element correlation, $\widetilde{\mathbf{c}}^{\nu_1\beta_{\nu_1}}\otimes\widetilde{\mathbf{c}}^{\nu_2\beta_{\nu_2}}$, with the coefficient $\mathbf{c}^{2,\nu_1\beta_{\nu_1},\nu_2\beta_{\nu_2}}$ of the basis vector $\ket{\nu_1\beta_{\nu_1},\nu_2\beta_{\nu_2}}$. Of course, any extra information of the system, e.g., symmetries of the kernels, can simplify the formal space.

With this interpretation, we can fully replace every term in \cref{eq:aoe} and get
\begin{equation}
    \mathbf{L}[\mathbf{u}] = \sum_{\mathrm{i}=0}^\infty\mathbf{F}_{\mathrm{i},\left[\nu_\mathrm{j}\beta_{\nu_\mathrm{j}}\right]}\odot\mathbf{c}^{\mathrm{i},\left[\nu_\mathrm{j}\beta_{\nu_\mathrm{j}}\right]} = \mathbf{0}
\end{equation}
as our OE, where 
\begin{equation}
    \ket{\mathbf{u}} = \sum_{\mathrm{i}=0}^\infty\mathbf{c}^{\mathrm{i}, \left[\nu_\mathrm{j}\beta_{\nu_\mathrm{j}}\right]}\ket{\left[\nu_\mathrm{j}\beta_{\nu_\mathrm{j}}\right]}_\mathrm{i}
\end{equation}
is the formal solution in $\Psi$ and 
\begin{equation}
    \mathbf{F}_{\mathrm{i}, \left[\nu_\mathrm{j}\beta_{\nu_\mathrm{j}}\right]} = \dfrac{1}{\mathrm{i}!}\int_{\Omega^\mathrm{i}}\mathbf{k}_{\mathrm{i}}^{\left[\eta_\mathrm{j}\right]}\odot\bigotimes_{\mathrm{j}=1}^\mathrm{i}\partial_{\mathbf{q}_\mathrm{j}}^{\eta_\mathrm{j}}b_{\nu_\mathrm{j}\beta_{\nu_\mathrm{j}}}\,\mathrm{d}\mathbf{q}_\mathrm{j}\,.
\end{equation}
Similarly, the weak formulation can be written as
\begin{equation}\label{eq:aug_weak}
    \sum_{\mathrm{i}=0}^\infty\mathbf{G}_{\mathrm{i},\beta,\left[\nu_\mathrm{j}\beta_{\nu_\mathrm{j}}\right]}\odot\mathbf{c}^{\mathrm{i}, \left[\nu_\mathrm{j}\beta_{\nu_\mathrm{j}}\right]} = \mathbf{0}
\end{equation}
where 
\begin{equation}
    \mathbf{G}_{\mathrm{i},\beta,\left[\nu_\mathrm{j}\beta_{\nu_\mathrm{j}}\right]}
    = \int_\Omega b^\ast_\beta\mathbf{F}_{\mathrm{i}, \left[\nu_\mathrm{j}\beta_{\nu_\mathrm{j}}\right]}\, \mathrm{d}\mathbf{q}\,.
\end{equation}
However, this finite system of algebraic equations has many more coefficients than \cref{eq:wf_fe_aoe}. Hence, extra conditions, which can take the form of regularization or full-blown equations, are required. 

\subsection{Tensor-Network Representation}
Due to the $\mathrm{i}$-linearity of $\mathbf{k}_{\mathrm{i}}^{\left[\eta_\mathrm{j}\right]}$, we have guage freedom to permute the pair indices $\left(\nu_\mathrm{j}, \beta_{\nu_\mathrm{j}}\right)$ of $\mathbf{G}_{\mathrm{i},\beta,\left[\nu_\mathrm{j}\beta_{\nu_\mathrm{j}}\right]}$ and $\mathbf{c}^{\mathrm{i}, \left[\nu_\mathrm{j}\beta_{\nu_\mathrm{j}}\right]}$. We chose to follow a canonical ordering such that for any $\mathrm{j}'>\mathrm{j}$, $\nu_{\mathrm{j}'}\geq\nu_\mathrm{j}$, i.e., in order of element labels. Furthermore, if $\nu_{\mathrm{j}'}=\nu_\mathrm{j}$, $\eta_{\mathrm{j}'}\geq\eta_\mathrm{j}$, i.e., in order of derivative order as well.

\begin{figure}[!htb]
    \centering
    \includegraphics[width=0.9\columnwidth]{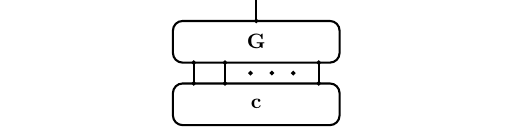}
    \caption{Tensor network representation of the left hand side of \cref{eq:TN_weak}. Unlike in \cref{fig:interact_kernel}, we omit the obvious tensor indices $\beta$ and $\left[\beta_{\nu\mathrm{j}}\right]$ for cleanliness.}
    \label{fig:residual}
\end{figure}

We can also tensorized the summation over $\mathrm{i}$ by introducing, for each element, empty basis function $b_{\nu0}$ into the local approximation space basis set in \cref{eq:local_approx_space}. With this addition and let $J$ be the largest correlation order contribution of the analytic OE from any element, we can map each basis vector in $\Psi$ to a new representation such that each element is pre-allocated with $J$ slots and the basis orders that appear in the basis vector are filled in the slots of their corresponding elements in the same order, e.g., if $J=3$ and $N=4$, the basis vector $\ket{13,14,21,22,33}$ is mapped to
\begin{equation*}
    \ket{(3,4,0),(1,2,0),(3,0,0),(0,0,0)}\,,
\end{equation*} 
i.e., 2 contributions from element 1, 2 contributions from element 2, 1 contribution from element 3, and no contribution from element 4. Hence, \cref{eq:aug_weak} can be replaced with
\begin{equation}\label{eq:TN_weak}
    \mathbf{G}_{\beta,\left[\beta_{\nu\mathrm{j}}\right]}\odot\mathbf{c}^{\left[\beta_{\nu\mathrm{j}}\right]} = \mathbf{0}
\end{equation}
with appropriate substitution of $\mathbf{G}_{\beta,\left[\beta_{\nu\mathrm{j}}\right]}$'s definition. The equation can be neatly represented in tensor network format as shown in \cref{fig:residual}.

Effectively, one can directly solve \cref{eq:TN_weak} and obtain the desired approximate solution. However, directly storing the full coefficient tensor $\mathbf{c}$ classically requires $O\left(B^{NJ}\right)$ coefficient blocks, where $B=\max\left(\{B_\mathrm{n}\}\right)$ which grows exponentially with the problem size. The full tensor-product space therefore cannot be introduced as a direct storage format, but as a formal correlated approximation space to be compressed by any TN ansatz.

\begin{figure}[!htb]
    \centering
    \includegraphics[width=\columnwidth]{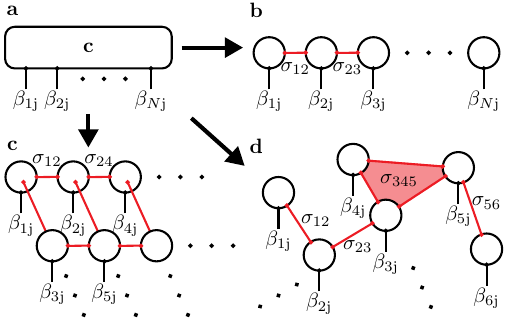}
    \caption{Tensor network ansatz of coefficient tensor with different locality assumptions: \textbf{(a)} no assumption, \textbf{(b)} matrix product state, \textbf{(c)} projected entangled pair state, and \textbf{(d)} higher-order  tensor network assumption where shaded regions represent index hyper-contraction. $\beta$'s and $\sigma$'s indicate tensor indices responsible for each bond and contraction, respectively.}
    \label{fig:tn_ansatz}
\end{figure}

The choice of TN ansatz is arbitrary in the perspective of our formalism, but, if one already has expectation for the solution to contain some correlation and/or entanglement behaviors, the network can be designed to support them. One crucial assumption is the locality of correlation: we expected the correlation to only appear between neighboring elements. \Cref{fig:tn_ansatz} illustrates representative TN ansatzes under different local correlation geometries.

Similar to \cref{fig:residual}, \cref{fig:tn_ansatz}a is the generic tensor representation of $\mathbf{c}$. \Cref{fig:tn_ansatz}b shows a matrix product state (MPS) which is appropriate for one-dimensional element orderings with nearest-neighbor correlations while \cref{fig:tn_ansatz}c shows a projected entangled pair state (PEPS) structure adapted for two-dimensional grid-like element connectivity. \Cref{fig:tn_ansatz}d illustrates a higher-order or hypergraph-type (HG) tensor network capable of representing multi-element correlations beyond pair-wise interactions. 

The bond dimensions (ranges of bond indices $\sigma$'s) control the amount of inter-element correlation that can be represented. If the maximum bond dimension in the TN is $S$, and the maximum bond degree (including hyper-bond) in the TN is $M$, then the storage requirement is 
\begin{equation}
    O(NB^JS^M)
\end{equation}
instead of $O(B^{NJ})$. Therefore, as long as the bond degree and correlation order are capped by a constant, the TN ansatz defines a polynomial-size variational manifold inside the exponentially large correlated coefficient space.

\subsection{Least-Squares Variational Formulation}
Given any formal coefficient tensor $\mathbf{c}$, the residual from weak formulation can be defined as
\begin{equation}
\mathbf{h}_{\beta}(\mathbf{c})
= \mathbf{G}_{\beta,\left[\beta_{\nu\mathrm{j}}\right]}\odot\mathbf{c}^{\left[\beta_{\nu\mathrm{j}}\right]}\,.
\end{equation}
We determine the correlated coefficients by minimizing the weighted residual norm
\begin{equation}
    J(\mathbf{c})
    = W^{\beta'\beta}\left(\mathbf{h}^\dagger_{\beta'}(\mathbf{c})\odot \mathbf{h}_{\beta}(\mathbf{c})\right)\,,
\end{equation}
where $W^{\beta'\beta}$ is a positive-definite metric on the residual space. The TN approximation is obtained by restricting $\mathbf{c}$ to a chosen TN ansatz variational manifold $\mathcal{T}_{\mathrm{TN}}$:
\begin{equation}
    \mathbf{c}_{\mathrm{TN}}
    =
    \arg\min_{\mathbf{c}\in\mathcal{T}_{\mathrm{TN}}}
    J(\mathbf{c}).
\end{equation}

Although this objective is quadratic and convex in the unrestricted correlated coefficient tensor $\mathbf{c}$, the TN parametrization makes the global optimization generally non-convex. However, when all tensors in TN except one are fixed, the residual becomes linear in the active tensor. Therefore, each single-site local update reduces the optimization to a convex weighted least-squares sub-problem which is part of the reason for the success of density matrix renormalization group (DMRG) algorithm for MPS ansatz.
\section{Specialization to Important Operator-Equation Classes}\label{sec:reduc}
The formulation above was stated for analytic operator equations in general form. In this section, we show how several common equation classes can be reduced to the same weak formulation structure instead of introducing different methods for different equations.

\subsection{Partial Differential Equations}
A partial differential equation (PDE) can be written as
\begin{equation}
    \mathbf{g}\left(\mathbf{q},\left[\partial^\eta_{\mathbf{q}}\mathbf{u}\right]\right)
    =
    \mathbf{0}\,.
\end{equation}
If $\mathbf{g}$ is analytic, then it can be expanded as a convergent series 
\begin{equation}\label{eq:reduce_PDE}
    \mathbf{g}\left(\mathbf{q},\left[\partial^\eta_{\mathbf{q}}\mathbf{u}\right]\right)
    = \sum_{\mathrm{i}}\dfrac{1}{\mathrm{i}!}\mathbf{a}_{\mathrm{i}}^{\left[\eta_\mathrm{j}\right]}\odot\bigotimes_{\mathrm{j}=1}^\mathrm{i}\partial_{\mathbf{q}}^{\eta_\mathrm{j}}\mathbf{u}\,,
\end{equation}
which is a special case of \cref{eq:aoe} when
\begin{equation}
    \mathbf{k}_{\mathrm{i}}^{\left[\eta_\mathrm{j}\right]} = \mathbf{a}_{\mathrm{i}}^{\left[\eta_\mathrm{j}\right]}\bigotimes_{\mathrm{j}=1}^{\mathrm{i}}\delta(\mathbf{q}_\mathrm{j}-\mathbf{q})\,.
\end{equation}

\subsection{Initial-Value Problems}\label{sec:reduce_IVP}
An initial-value problem (IVP) is simply a variation of PDE with time as one of the domain dimensions. Hence, it is of the form
\begin{equation}
    \mathbf{g}\left(t, \mathbf{q}, \left[\partial^\eta_{\mathbf{q}}\partial^\zeta_{t}\mathbf{u}\right]\right)
    =
    \mathbf{0}\,.
\end{equation}
If $\mathbf{g}$ is analytic, 
\begin{multline}
    \mathbf{g}\left(t, \mathbf{q}, \left[\partial^\eta_{\mathbf{q}}\partial^\zeta_{t}\mathbf{u}\right]\right)\\
    = \sum_{\mathrm{i}}\dfrac{1}{\mathrm{i}!}\mathbf{a}_{\mathrm{i}}^{\left[\zeta_\mathrm{k}\right], \left[\eta_\mathrm{j}\right]} \odot\bigotimes_{\mathrm{j}=1}^\mathrm{i}\bigotimes_{\mathrm{k}=\mathrm{j}+1}^\mathrm{i}\partial_{\mathbf{q}}^{\eta_\mathrm{j}}\partial_{t}^{\zeta_\mathrm{k}}\mathbf{u}\,.
\end{multline}

For IVP, given $\mathbf{u}'$ and some of its time derivatives at time $t$, the problem is to determine $\mathbf{u}$ and those time derivatives at $t+\Delta t$. There are many schemes one can use with varying degree of convergence criteria and computational costs. We will use a simple implicit Euler's scheme as an example that replaces the analytic form of the IVP with
\begin{multline}
    \mathbf{0} = \sum_{\mathrm{i}}\dfrac{1}{\mathrm{i}!}\mathbf{a}_{\mathrm{i}}^{\left[\zeta_\mathrm{k}\right],\left[\eta_\mathrm{j}\right]}\odot\bigotimes_{\mathrm{j}=1}^\mathrm{i}\bigotimes_{\mathrm{k}=\mathrm{j}+1}^\mathrm{i}\\
    \partial_{\mathbf{q}}^{\eta_\mathrm{j}}
    \begin{cases}
        \dfrac{1}{\Delta t}\left(\partial_{t}^{\zeta_\mathrm{k}-1}\mathbf{u}-\partial_{t}^{\zeta_\mathrm{k}-1}\mathbf{u}'\right), &\zeta_\mathrm{k}\neq 0\\
        \mathbf{u}, &\zeta_\mathrm{k} = 0
    \end{cases}\,.
\end{multline}
With some term arrangement, the equation can be written as a PDE on $\mathbf{v}=\left[\partial_{t}^{\zeta}\mathbf{u}\right]$, a special case of \cref{eq:aoe}.

\subsection{Integro-Differential Equations}
An integro-differential equation (IDE) is a generalization of PDE to include non-local interaction mostly in the form of integration kernel as
\begin{equation}
    \mathbf{g}\left(\mathbf{q},\left[\partial^\eta_{\mathbf{q}}\mathbf{u}\right], \left[\int_{\Omega^\mathrm{i}}\mathbf{a}_\mathrm{i}\left(\mathbf{q}, \left[\mathbf{q}_\mathrm{j}\right], \left[\partial^{\eta}_{\mathbf{q}_\mathrm{j}}\mathbf{u}_\mathrm{j}\right]\right)\,\prod_{\mathrm{j}=1}^\mathrm{i}\mathrm{d}\mathbf{q}_\mathrm{j}\right]\right)
    =
    \mathbf{0}\,.
\end{equation}
If $\mathbf{g}$, and all $\mathbf{a}_\mathrm{i}$'s are analytic, then it
can be expanded as a special case of \cref{eq:aoe} with multi-linear interaction kernel containing some Dirac delta functions for the parts that are not integrations in the IDE.

\subsection{Memory and Delay Equations}
A memory-delay equation (MDE) is an IDE with non-local interaction in temporal dimension, i.e.,
\begin{multline}
    \mathbf{0} = \mathbf{g}\Biggl(t, \mathbf{q},\left[\partial^\eta_{\mathbf{q}}\partial^\zeta_{t}\mathbf{u}\right],\\
    \left.\left[\int_{\tau^\mathrm{i}}\mathbf{a}_\mathrm{i}\left(t, \left[t_\mathrm{j}\right], \mathbf{q}, \left[\partial^{\eta}_{\mathbf{q}}\partial^{\zeta}_{t_\mathrm{j}}\mathbf{u}_\mathrm{j}\right]\right)\,\prod_{\mathrm{j}=1}^\mathrm{i}\mathrm{d}t_\mathrm{j}\right]\right)\,.
\end{multline}
Similar to the IDE, if $\mathbf{g}$, and all $\mathbf{a}_\mathrm{i}$'s are analytic, then it
can be expanded as a special case of \cref{eq:aoe}.
\section{Demonstration}\label{sec:demo}
The formalism developed above applies to a broad class of analytic OEs across different geometries and TN ansatzes. To illustrate the framework concretely, we specialize to one-dimensional (1D) problems which can be naturally discretized into an MPS variational problems since there is a well-behave solving algorithm (DMRG) for them. We demonstrate the approach on IVP and examine the convergence and accuracy of the resulting TN solutions.

\subsection{Specialization}
As a continuation from \cref{sec:reduce_IVP}, the ending note of that sub-section implies that we only need to consider the IVPs with only first order time derivatives. Also, instead of $\mathbf{u}'$ and $\mathbf{u}$, we will use $\mathbf{u}^{(\mathrm{t}-1)}$ and $\mathbf{u}^{(\mathrm{t})}$ notations where labels in parentheses indicate time step. Hence, for 1D, we are ultimately solving for $\mathbf{u}^{(\mathrm{t})}$ in
\begin{multline}
    0 = \sum_{\mathrm{i}}\dfrac{1}{\mathrm{i}!}a_{\mathrm{i}}^{\left[\zeta_\mathrm{k}\right],\left[\eta_\mathrm{j}\right]}\prod_{\mathrm{j}=1}^\mathrm{i}\prod_{\mathrm{k}=\mathrm{j}+1}^\mathrm{i} \\
    \partial_{q}^{\eta_\mathrm{j}}
    \begin{cases}
        \dfrac{1}{\Delta t}\left(u^{(\mathrm{t})}-u^{(\mathrm{t}-1)}\right), &\zeta_\mathrm{k} = 1\\
        u^{(\mathrm{t})}, &\zeta_\mathrm{k} = 0
    \end{cases}\,.
\end{multline}
which is a PDE but not in the same format as \cref{eq:reduce_PDE}: the existence of $u^{(\mathrm{t}-1)}$ which can be treated as part of another term with lower correlation order of $u^{(\mathrm{t})}$.

It is straight forward to simply use the weak formalism directly on the equation, but it means that we need to calculate $\mathbf{G}$ from scratch at every time step since $u^{(\mathrm{t}-1)}$ got updated. Instead, we can leverage the deterministic structure of the IVP to pre-compute some parts. Different classes of OEs will be different but for IVE and, to extension, PDE the spatial derivatives and multiplications can be pre-formuated and approximated into concise weak forms. 

As spatial derivatives and multiplications are local operators, with well-distributed basis function, any of their actions on $\Phi_\mathrm{n}$ can be approximated with a member of $\Phi_\mathrm{n}$. In other words, we can approximate the spatial derivatives and multiplications with some tensor contractions. 

\begin{figure}[!htb]
    \centering
    \includegraphics[width=0.9\columnwidth]{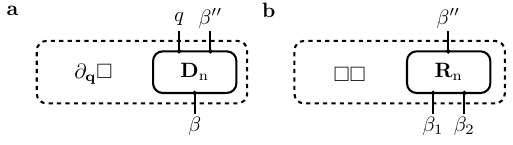}
    \caption{Tensor circuit representation of \textbf{(a)} spatial derivative and \textbf{(b)} multiplication act locally at element $\mathrm{n}$.}
    \label{fig:dif_mul_tn}
\end{figure}

Consider an arbitrary spatial derivative of a basis function $b_{\mathrm{n},\beta}$ which can be approximated as
\begin{equation}
    \partial_{\mathbf{q}}b_{\mathrm{n},\beta} \approx \mathbf{m}^{\beta''}b_{\mathrm{n},\beta''}\,,
\end{equation}
which, by the weak formulation, gives
\begin{equation}
    \int_\Omega b^\ast_{\mathrm{n},\beta'}\partial_{\mathbf{q}}b_{\mathrm{n},\beta}\,\mathrm{d}q \approx \mathbf{m}^{\beta''}\int_\Omega b^\ast_{\mathrm{n},\beta'}b_{\mathrm{n},\beta''}\,\mathrm{d}q\,.
\end{equation}
We then define
\begin{align}
    \mathbf{P}_{\mathrm{n},\beta'\beta} &= \int_\Omega b^\ast_{\mathrm{n},\beta'}\partial_{\mathbf{q}}b_{\mathrm{n},\beta}\,\mathrm{d}q\, ,\\
    Q_{\mathrm{n},\beta'\beta''} &= \int_\Omega b^\ast_{\mathrm{n},\beta'}b_{\mathrm{n},\beta''}\,\mathrm{d}q\, ,
\end{align}
which can be pre-computed given the basis functions. If $Q^{-1}_{\mathrm{n},\beta''\beta'}$ is the matrix inverse of $Q_{\mathrm{n},\beta'\beta''}$. Then, $\partial_{\mathbf{q}}$ on the element $\mathrm{n}$ can be approximated with
\begin{equation}
     \mathbf{D}_{\mathrm{n},\beta''\beta}\approx Q^{-1}_{\mathrm{n},\beta''\beta'}\mathbf{P}_{\mathrm{n},\beta'\beta}\,.
\end{equation}
We can also visualize it in the tensor circuit form as shown in \cref{fig:dif_mul_tn}a. Note that the label $\mathrm{n}$ is needed since each element has different sets of besis functions. Also, the index $q$ of the derivative tensor indicate the axis of differentiation.

Similarly, consider an arbitrary spatial multiplication of a pair of basis functions $b_{\mathrm{n},\beta_1}$ and $b_{\mathrm{n},\beta_2}$ which can be approximated as
\begin{equation}
    b_{\mathrm{n},\beta_1}b_{\mathrm{n},\beta_2} \approx m^{\beta''}b_{\mathrm{n},\beta''}\,.
\end{equation}
If we define,
\begin{align}
    P'_{\mathrm{n},\beta'\beta_1\beta_2} &= \int_\Omega b^\ast_{\mathrm{n}\,\beta'}b_{\mathrm{n},\beta_1}b_{\mathrm{n},\beta_2}\,\mathrm{d}q\,,\\
    Q'_{\mathrm{n},\beta'\beta''} &= \int_\Omega b^\ast_{\mathrm{n},\beta'}b_{\mathrm{n},\beta''}\,\mathrm{d}q\,,
\end{align}
then, any multiplication in this elemental basis can be approximated with tensor product reduction
\begin{equation}
    R_{\mathrm{n},\beta''\beta_1\beta_2} = Q^{\prime -1}_{\mathrm{n},\beta''\beta'}P'_{\mathrm{n},\beta'\beta_1\beta_2}\,,
\end{equation}
as depicted in \cref{fig:dif_mul_tn}b. 

\begin{figure}[!htb]
    \centering
    \includegraphics[width=0.9\columnwidth]{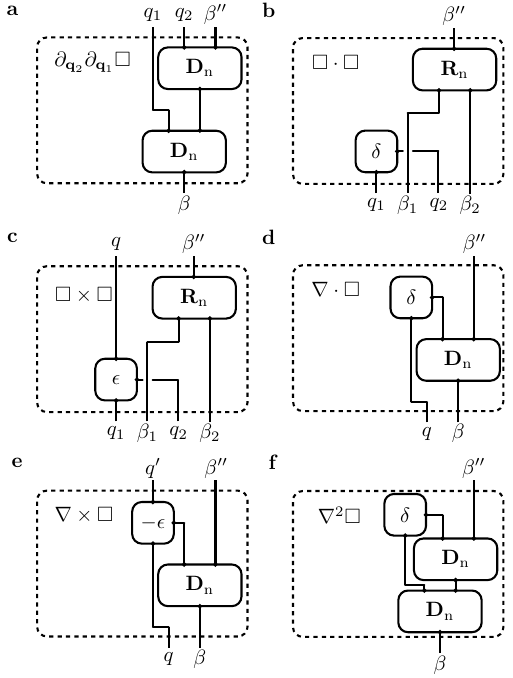}
    \caption{Tensor circuit representation of more complicated operators can be constructed from the two basic building block types (spatial derivative and multiplication) together with well known constant tensors like Kronecker delta ($\delta$) and Levi-Civita ($\epsilon$). \textbf{(a)} Second spatial derivative or Hessian can be achieve with stacking two $\mathbf{D}$'s with both spatial indices free. \textbf{(b)} Dot and \textbf{(c)} cross products are multiplications while \textbf{(d)} Divergence and \textbf{(e)} curl are first spatial derivatives with additional spatial index summation by $\delta$ and $\epsilon$. Lastly, \textbf{(f)} Laplacian is second spatial derivatives with additional spatial index summation by $\delta$. The subscript $\mathrm{n}$ make it explicitly clear that the derivatives and multiplications act locally at element $\mathrm{n}$.}
    \label{fig:complex_tn}
\end{figure}

With these components, we can construct other more complicated IVP operators, e.g., second derivative, dot product, cross product, divergence, curl, and Laplacian, as illustrated in \cref{fig:complex_tn}. Of course, directly applying the weak formalism for each operator would give a better approximation at the cost of additional overhead calculations per time step. 

\subsection{Example}
We'll apply the procedure of converting IVPs to TN with the following non-linear diffusion equations
\begin{equation}
    \partial_tu = u^\mathrm{p}\partial^2_xu
\end{equation}
with $u = u(t, x)$, $\mathrm{p}\in\{0, 1, 2\}$, $\Omega = [0, 1]$, $u(0, x)$ is a modified Gaussian centered at $x = 0.5$ with $\partial_xu(t, 0)=\partial_xu(t, 1)=0$, i.e., close boundaries. We can rewrite the IVP in our format as
\begin{equation}\label{eq:nlIVP}
    0 = u^{(\mathrm{t})\mathrm{p}}\partial^2_xu^{(\mathrm{t})} - \dfrac{1}{\Delta t}(u^{(\mathrm{t})}-u^{(\mathrm{t}-1)})\,.
\end{equation}

\begin{figure}[!htb]
    \centering
    \includegraphics[width=\columnwidth]{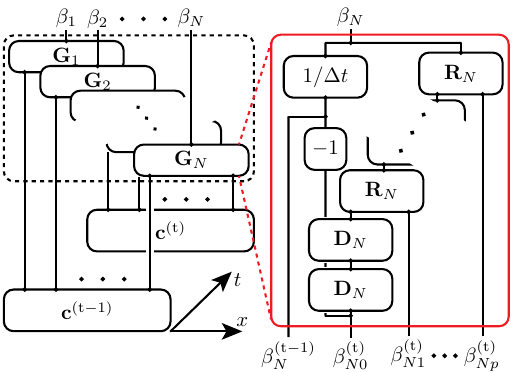}
    \caption{Tensor circuit representation of residual of non-linear diffusion IVP. Since IVPs only contain local interactions, the circuit elements ($\mathbf{G}_\mathrm{n}$) only connect branches at the same elements. The circuit in red box illustrates the detail of $\mathbf{G}_\mathrm{N}$. Since the first term in \cref{eq:nlIVP} is $u^{(\mathrm{t})\mathrm{p}}\partial^2_xu^{(\mathrm{t})}$, there are 2 derivatives tensor on $\mathbf{c}^{(\mathrm{t})}$ and $p$ multiplications between $p$ $\mathbf{c}^{(\mathrm{t})}$'s and output branch from the second derivative of $\mathbf{c}^{(\mathrm{t})}$ in the red box. The boxes containing $-1$ and $1/\Delta t$ implies multiplication by those values and branch merging implies addition of merged branches which corresponding to the second term in \cref{eq:nlIVP}.}
    \label{fig:diffuse_tc}
\end{figure}

The corresponding tensor circuit of this IVP is as shown in \cref{fig:diffuse_tc}.

\begin{figure*}[!htb]
    \centering
    \includegraphics[width=0.9\textwidth]{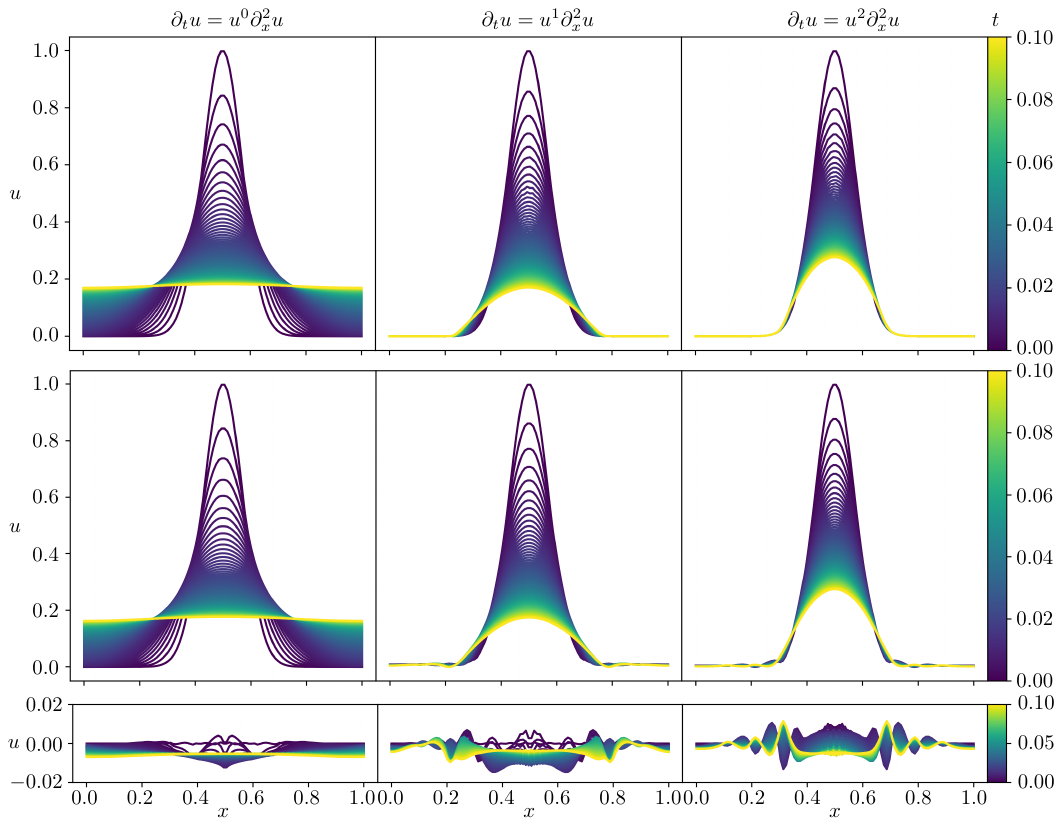}
    \caption{Numerical simulation results for a non-linear diffusion IVPs (\cref{eq:nlIVP}) with different $p$'s: $0$ (left column), $1$ (middle column), $2$ (right column). Top row shows result from explicit Runge–Kutta method of order 5(4)\cite{DORMAND198019} while middle row are the results from our TN variational method with $N=11$, $B=4$, $S=1$ with their differences is shown in the bottom row..}
    \label{fig:diffuse}
\end{figure*}

Then, we solve the non-linear diffusion IVP for each exponent $\mathrm{p}$ using $N$ uniform finite elements on $\Omega$, each equipped with $B$ local Lagrange polynomial basis functions, and an MPS ansatz with bond dimension $S$ for the coefficient tensor. The time evolution is advanced using the implicit Euler scheme described with time step $\Delta t$. At each time step, the DMRG sweep algorithm minimizes the weighted residual $J(\mathbf{c})$ over the MPS variational manifold $\mathcal{T}_{\mathrm{MPS}}$.

To ensure the continuity and smoothness of the results, we also add an additional circuit to measure the residual of multiple derivative orders at each boundary between elements as well as the residual of Neumann boundary conditions at the outer edges of the domain to environment. As these continuities only involve at most 2 adjacent elements, the circuit blocks either only connect adjacent tensor index branches (element boundaries) or only connect to one branch (domain edges).

\Cref{fig:diffuse} summarize the time-evolution results across all three cases for  explicit Runge-Kutta IVP simulation (top row) and TN variational (middle row) method. For $\mathrm{p}=0$ (left column), the equation reduces to the standard diffusion equation $\partial_t u = \partial_x^2 u$, which admits a smooth, rapidly decaying solution. For $\mathrm{p}=1$ (middle column), the non-linear coupling $u\partial_x^2 u$ introduces amplitude-dependent diffusion, slowing the spread of $u$ in regions of low amplitude. Finally, for $\mathrm{p}=2$ (right column), the stronger non-linearity further suppresses diffusion and produces sharper features. In all cases, the TN variational method with $N=10$, $B=4$, $S=1$ produce a reliable results with the difference of less than two percent of the original amplitude of the modified Gaussian as shown in the bottom row of \cref{fig:diffuse}. We also see that the Neumann boundary conditions at both ends of the domain as well as the continuity and smoothness of the results are also respected to within numerical precision confirming that the weak formulation correctly encodes the boundary constraints.

\subsection{Convergence and Accuracy}
\begin{figure}[!htb]
    \centering
    \includegraphics[width=0.9\columnwidth]{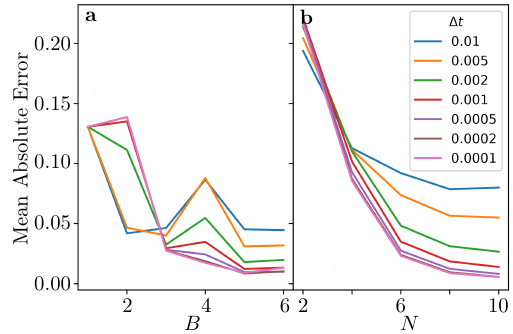}
    \caption{Mean absolute error trend of TN variational method for solving diffusion IVP with different \textbf{(a)} basis order $B$ and \textbf{(b)} number of elements $N$ at different time step size $\Delta t$. For each fixed $\Delta t$, the trends exhibit convergence up to a level after which the better temporal resolution is required for a more accurate result.}
    \label{fig:convergence}
\end{figure}

To assess the convergence properties of the TN solver, we study the behavior of the solution error as functions of Lagrange polynomial basis order $B$, and the number of finite elements $N$ at different time step size $\Delta t$. Since both $B$ and $N$ are, in some way, representation of spatial resolution, we should expect the behavior of the trend to be similar to standard finite difference method which exhibit the stability condition where the convergence cannot be improved purely by increasing spatial resolution while the temporal resolution is fixed. This is clearly shown in \cref{fig:convergence}.

Together, these results demonstrate that the TN framework introduced in this paper can produce accurate, converged solutions to non-linear IVPs with the computational scaling controlled by the bond dimension rather than the full exponential size of the coefficient space. While the examples here are restricted to 1D MPS/DMRG for clarity, the formalism is directly extensible to higher-dimensional problems via PEPS or tree tensor network ansatzes.
\section{Conclusion}\label{sec:conc}
We have introduced a tensor-network formulation for finite-element discretizations of analytic operator equations. Starting from the Galerkin weak formulation, analytic operator equations generate a hierarchy of multilinear residual tensors determined by the operator kernels, basis functions, and element geometry. By lifting the conventional finite-element coefficient space into an augmented tensor-product space over elements, the discretized problem can be expressed as a variational optimization over correlated coefficient states. Tensor networks then provide a structured parametrization of this enlarged space, allowing inter-element correlations to be represented without explicitly storing the exponentially large full coefficient tensor.
The key distinction of this framework is that tensor networks are used not only as compression tools for already-discretized solutions, but as variational representations of the finite-element coefficient space itself. This perspective establishes a common formulation for a broad class of analytic operator equations, including partial differential equations, initial-value problems, integro-differential equations, and memory or delay equations. Different equation classes enter through different interaction kernels and residual tensors, while the variational tensor-network structure remains the same.
For initial-value problems, we showed that implicit time discretization reduces the evolution to a sequence of stationary variational problems. Spatial derivative and multiplication operators can be precomputed at the element level and reused across time steps, reducing repeated assembly costs. The one-dimensional nonlinear diffusion examples demonstrate that an MPS-based formulation can reproduce standard numerical solutions with controlled accuracy while respecting continuity and boundary constraints.
The present results should be viewed as a proof of principle rather than a guarantee of universal computational advantage. The efficiency of the method depends on whether the solution admits a compact tensor-network representation. Highly discontinuous, chaotic, or strongly correlated solutions may require large bond dimensions, and accuracy remains sensitive to both local basis quality and the chosen network geometry. Future work should investigate adaptive bond-dimension control, higher-dimensional tensor-network ansatzes such as PEPS and tree networks, rigorous convergence estimates, and comparisons with conventional finite-element solvers on large-scale benchmark problems.
Overall, the proposed framework provides a direct bridge between finite-element approximation theory, multi-linear operator discretization, and tensor-network variational optimization. It offers a flexible foundation for classical, quantum-inspired, and potentially quantum-assisted solvers for structured operator equations.
\section*{Acknowledgments}\label{sec:ack}
The authors thank Alan Edelman and Nuno Loureiro for insightful discussions. We dedicate this work to  Loureiro's memory and note with sadness that he passed away during the course of this work. This work was supported by the National Science Foundation (NSF) Convergence Accelerator Award No. 2345084, the U.S. Department of Energy (DOE), Basic Energy Sciences (BES), Award No. DE-SC0020148, MIT Energy Initiative, and support from R. Wachnik.
\appendix

\bibliographystyle{unsrt}
\bibliography{references}
 
\end{document}